\documentclass[12pt]{amsart}

\usepackage{amsmath}
\usepackage{amssymb}
\usepackage{amscd}
\usepackage{amsthm}
\usepackage[dvips]{graphics}

\def\no{\noindent}

\theoremstyle{plain}
\newtheorem{theorem}[equation]{Theorem}

\newtheorem{proposition}[equation]{Proposition}

\newtheorem{lemma}[equation]{Lemma}
\newtheorem{corollary}[equation]{Corollary}

\theoremstyle{remark}

\theoremstyle{definition}

\newcommand{\ra}{\rightarrow}

\newcommand{\N}{\mathbb N}

\newcommand{\R}{\mathbb R}

\renewcommand{\H}{\mathbb H}

\renewcommand{\S}{\mathbb S}

\renewcommand{\:}{\colon}

\newcommand{\sph}{\S^1}
\newcommand{\sub}{\subset}

\newcommand{\al}{\alpha}

\def\eps{\epsilon}

\def\la{\lambda}

\def\ra{\rightarrow}

\def\geo{\partial_{\infty}}

\def\defeq{:=}

\newcommand{\qs}{quasisymmetric}

\begin{document}

\title{Quasi-hyperbolic planes in hyperbolic groups}

\author{Mario Bonk}
\address{Department of Mathematics, University of Michigan,
Ann Arbor, MI, 48109-1109}
\email{mbonk@umich.edu, bkleiner@umich.edu}
\author{Bruce Kleiner}

\thanks{M.B.\ was  supported by NSF grant DMS-0200566. 
B.K.\ was supported by NSF grants  DMS-9972047 and DMS-0204506.}

\date{January~19, 2003}

\begin{abstract}
The hyperbolic plane
$\H^2$ admits a quasi-isometric embedding into
every hyperbolic group which is not virtually free. \end{abstract}

\maketitle

\noindent 
The purpose of this note is to prove the following theorem which 
answers a question posed by P.~Papasoglu:
\begin{theorem}
\label{bigtheorem}
The hyperbolic plane
$\H^2$ admits a quasi-isometric embedding into
a  hyperbolic group if and only if the group is not
virtually free.  
\end{theorem}

\noindent 
A map $f\: X \ra Y$ between two metric spaces $(X,d_X)$ and $(Y, d_Y)$
 is called   a {\em quasi-isometric embedding} if there exist
constants $\la \ge 1$ and $K \ge 0$ such that
$$ \frac 1\la d_X(x,y)-K \le d_Y(f(x), f(y))\le \la d_X(x,y)+K$$
for all $x,y\in X$.
A group is  {\em virtually free} if it contains a free 
subgroup of finite index.    
We refer to \cite{ghysdela} for the
definition of hyperbolic groups and related
concepts from the theory of Gromov hyperbolic spaces.
Every Gromov hyperbolic space $X$ 
 has a boundary $\partial_\infty X$ 
which carries a class of canonical {\em visual   metrics}. 
These metrics are bi-Lipschitz equivalent to distance functions 
of the form 
$$ d_{w,\eps}(a,b)= \exp(-\eps (a,b)_w), \quad a,b\in \partial_\infty X, $$
where $w\in X$ is a base point, $\eps>0$ is sufficiently small, and 
$(a,b)_w$ denotes the Gromov product of the points $ a$ and $b$
with respect to $w$ (cf.\ \cite[Ch.~7]{ghysdela}).

\begin{corollary} \label{cor}
The boundary of  a hyperbolic group
(equipped with any visual metric) contains a quasi-circle
if and only if the group is not virtually free.
\end{corollary}
   
\no 
By definition a {\em quasi-circle} is a metric circle which admits a
quasisymmetric parametrization by the unit  circle $\sph \sub \R^2$
(see \cite{hein} for the definition and  basic
facts about quasisymmetric maps).   Since the boundary of a virtually
free group is totally disconnected, the ``only if'' part of the 
corollary is obvious.

One of the main ingredients in the proof of the theorem is a result 
by Tukia \cite{tukia} which insures the existence of quasi-arcs with given
end-points inside certain  subsets of $\R^n$ (a
{\em quasi-arc} is a quasisymmetric image of the interval
$[0,1]$). The authors 
 would like to thank
Juha Heinonen for drawing their  attention to  Tukia's
paper, which allowed them  to substantially shorten
the proof of the next    proposition.

To state this proposition, we need one more  definition. 
A metric space $Z$ is {\em linearly connected} if
there exists  a constant $L$ such that for all $x,y\in Z$
there is a connected subset $S\subset Z$ of diameter
at most $Ld(x,y)$ containing $\{x,y\}$.

\begin{proposition}
\label{qarcsexist}
If $X$  is a complete, doubling,  and
linearly  connected metric space, then any two
distinct  points in $X$ are the endpoints of a quasi-arc. 
\end{proposition}
\proof
Let $d$ denote metric on $X$, and pick $\al\in (0,1)$.
Since $X$ is doubling, 
there exists $n\in \N$ such that the ``$\al$-snowflaked'' metric
space $(X,d^\al)$ can be embedded  into $\R^n$ (equipped with the usual metric)
 by a bi-Lipschitz mapping (this follows from Assouad's Embedding
Theorem \cite[2.6.~Prop.]{assouad}; see \cite[Thm.~12.2]{hein}
for the version of this theorem used here).
Let $Z$ denote the image of such an embedding.
Then $Z$ is complete and linearly
 connected, since $X$ has these properties. 
Hence 
 any two distinct  points in $Z$ are the endpoints of a quasi-arc in $Z$
(up to terminology this  is  \cite[Thm 1A]{tukia};  see the introduction of
 \cite{tukia} for a discussion). 
Since quasi-arcs in $Z$ pull back to quasi-arcs in $X$,  the result 
follows.  \qed

\begin{proposition}
\label{linloccon}
If  $G$ is  $1$-ended hyperbolic group, 
then 
$\geo G$ equipped with any visual metric  $d$
  is compact, doubling, connected, and linearly  
connected.
\end{proposition}
\proof It is easy  to show that $\geo G$ is compact
\cite[p.~123, 9.~Prop.] {ghysdela} and doubling \cite[Sect.~9]{bonkschramm}.
Since the group $G$ is  $1$-ended, 
its boundary $\geo G$ is connected. 

It remains to  prove linear  connectedness (note that
this a stronger quantitative
version of  local connectedness which  was established  in this context 
in  \cite[Prop. 3.3]{bestmess}).
  Given two points $x$ and $y$
in a metric space $(Z,d)$, and $\la>0$, a {\em $\la$-chain from 
$x$ to $y$}  is a sequence of points $x=z_1,\dots,z_k=y$
such that $d(z_i,z_{i+1})\leq \la$ for all $1\leq i<k$.
The {\em length} of a $\la$-chain is the number of points
in the chain.

\begin{lemma}
There is a  number $N\in \N$ such that for all $x,\,y\in \geo G$ there is a
$\frac{1}{2}d(x,y)$-chain of length at most $N$ from $x$ to $y$.
\end{lemma}
\proof
If not, there are sequences $\{x_k\},\,\{y_k\}\subset\geo G$ such that the
shortest $\frac{1}{2}d(x_j,y_j)$-chain from $x_j$ to $y_j$
has length $j$.  The boundary $\geo G$ is compact and connected,
so clearly $r_j\defeq d(x_j,y_j)\ra 0$  as $j\ra\infty$. 
In view of the doubling property,
the sequence $(\geo G,\frac1{r_j}d,x_j)$
of pointed metric spaces subconverges to  a limit  $(W,d_W,x_\infty)$
with respect to  pointed Gromov-Hausdorff convergence
\cite[Thm.~8.1.10]{BBI}. We can then find a 
point $y_\infty\in W$ such that $d_W(x_\infty,y_\infty)=1$
and there is no $\la$-chain from $x_\infty$ to $y_\infty$
for any $\la<\frac{1}{2}$.   This implies that 
$W$ is not connected. By \cite[Lemma 5.2]{quasimobius},
the limit space $W$ is homeomorphic to $\geo G\setminus\{z\}$ for 
some $z\in\geo G$, and so $z$ is a ``global cut point" of $\geo G$. 

On the other hand, it is a well-known (and deep) fact 
 if $\geo G$ is connected, then $\geo G$ has no 
global cut points (see  \cite{swarup}, \cite[Thm.~9.3]{bowditch5},
\cite[Cor.~0.3]{bowditch8}). 
This is a contradiction.
\qed 

\medskip
\no
Now suppose $x$ and $y$ are arbitrary 
 points in $\geo G$.  By the lemma we can find
a  $\frac{1}{2}d(x,y)$-chain $S_1=\{z_1, \dots, z_k\}$ which joins 
$x$ to $y$ and has length $k\le N$.
Now define $S_2$ by adding, for  each $1\leq i<k$, 
the points in a $\frac{1}{2}d(z_i,z_{i+1})$-chain
joining $z_i$ to $z_{i+1}$.  Repeating  this process inductively, 
we  obtain 
a nested sequence of sets $S_1\subset\ldots\subset S_j\subset\ldots$.
The closure $S$ of the union $\bigcup_j S_j$ will be a connected set 
containing $x$ and $y$ whose diameter does not exceed $Ld(x,y)$, 
where $L$ is a constant independent of $x$ and $y$.
This shows that $\geo G$ is linearly connected. 
\qed

\medskip
\no
{\em The proofs of Theorem \ref{bigtheorem} and Corollary~\ref{cor}.}
We first assume that $G$ a hyperbolic group which is not virtually free,
and prove that there is a quasi-isometric embeddding $\H^2\ra G$
and a quasi-circle in $\geo G$.
Every hyperbolic group is
finitely presentable \cite[p.~76, 17.~Prop.]{ghysdela}.
Hence 
there is a finite graph of groups
decomposition of  $G$ where all edge groups are finite,
and all vertex groups have at most one end \cite[Theorem 6.2.14]{zieschang}.
Since  $G$ is not virtually free, 
one of the vertex groups $G_0$ is $1$-ended
\cite[Theorem 6.2.12]{zieschang}. The group $G_0$ 
is quasi-isometrically embedded in $G$, since this is true for every 
vertex group in a  graph of groups
decomposition with finite  edge groups
\cite[Rem.~3.6]{kapovich}. 
This 
implies  that $G_0$
is also a hyperbolic group. 
So
without loss of generality we may assume that $G$ itself
is $1$-ended.

Let $\geo G$
denote the boundary of $G$ equipped with a visual metric.
By Proposition~\ref{linloccon}, the hypotheses of Proposition
\ref{qarcsexist} are satisfied for $\geo G$.  Hence there is a \qs\ map 
$[0,1]\ra\geo G$.  Since $[0,1]$ is \qs ally homeomorphic
to the boundary of 
a hyperbolic half-plane $\H^2_+\subset \H^2$, we
 conclude that there is a quasi-isometric 
embedding $\H^2_+\ra G$ (see the proof of  Prop.~4.2
in \cite{Pau}, for example).  In particular,  one can 
quasi-isometrically embed arbitrarily large balls
$B\subset \H^2$ into $G$ with uniform constants for 
the  quasi-isometric embeddings. 
 By pre-composing with isometries
in $\H^2$, post-composing with left translations in the
group $G$, and applying a  compactness argument based on the
 Arzel\`a-Ascoli Theorem,
we can obtain a quasi-isometric embedding $\H^2\ra G$ 
as a limit.  A quasi-isometric embedding
of a Gromov hyperbolic space $X$ into a Gromov hyperbolic space 
$Y$ induces a quasisymmetric embedding of $\geo X$ into $\geo Y$ (see 
\cite[Thm.~6.5]{bonkschramm}, where this is essentially proved);
since $\geo \H^2$ is quasisymmetrically equivalent to $\sph$, we
deduce that the boundary
$\geo G$ contains a quasi-circle.  

Now suppose $G$ is virtually free.  It follows that $\geo G$
is totally disconnected, and therefore cannot contain a quasi-circle.
This then implies that there is no quasi-isometric embedding
$\H^2\ra G$.  

This completes the proofs of the theorem and corollary.
\qed

\medskip 
\no
{\em Remarks.} There are various open questions that are related 
to our theorem. For example,  
Papasoglu has  asked if every one-ended finitely
presented group $G$ contains a quasi-plane---the image of a
uniform embedding $P\ra G$ where $P$ is a complete
Riemannian plane of bounded geometry.
A problem due to Gromov is whether  every $1$-ended
hyperbolic group $G$ is the target of a homomorphism $\phi:S\ra G$
where $S$ is a surface group and $\phi$ does not factor
through a free group.
 

\bibliography{refs}
\bibliographystyle{siam}
\addcontentsline{toc}{subsection}{References}

\end{document}